\documentstyle{article}
\newcommand{\la}{\lambda}
\newcommand{\vsp}{\medskip}
\newcommand{\be}{\begin{equation}}
\newcommand{\ee}{\end{equation}}
\newcommand{\bea}{\begin{eqnarray}}
\newcommand{\eea}{\end{eqnarray}}
\newcommand{\nn}{\nonumber}
\newcommand{\Pe}{{\cal P}}
\newcommand{\ra}{\rightarrow}
\newcommand{\I}{(i_1|\cdots |i_r)}
\newcommand{\J}{(j_1|\cdots |j_s)}
\newcommand{\K}{(k_1|\cdots |k_t)}
\newcommand{\ot}{\otimes}
\newcommand{\op}{\oplus}
\newcommand{\tip}{\tilde{\pi}}
\newcommand{\ar}{\{ a^{(1)}|\cdots |a^{(r)}\}}
\newcommand{\as}{\{ a^{(1)}|\cdots |a^{(s)}\}}
\newcommand{\att}{\{ a^{(1)}|\cdots |a^{(t)}\}}
\newcommand{\zip}{z_{\tip}(\tip)}
\newcommand{\defi}{\stackrel{\rm def}{=}}
\newcommand{\CA}{C^{\bullet}(A)}
\newcommand{\Ai}{A_{\infty}}
\newcommand{\Li}{L_{\infty}}
\newcommand{\Gi}{G_{\infty}}
\newcommand{\KM}{K_{\bullet}\underline{{\cal M}}}
\newcommand{\usi}{u_1+\cdots +u_s=i-1}
\newcommand{\vtj}{v_1+\cdots +v_t=j_{\alpha}+u_{\alpha}-1}
\newcommand{\wtj}{w_1+\cdots +w_t=j_{\alpha}-1}
\newcommand{\zsi}{z_1+\cdots +z_{s+t-1}=i-1}
\newcommand{\tmtm}{\{\tilde{m}\}\{\tilde{m}\} =0}

\newcommand{\De}{\Delta}
\newcommand{\abc}{\{ a|b,c|d,e,f\}}
\newcommand{\Peb}{\bar{\Pe}}
\newcommand{\Be}{{\cal B}}
\newcommand{\Beb}{\bar{\Be}}
\newcommand{\CPA}{C_{\Pe}^{\bullet}(A)}
\newcommand{\CAA}{C^{\bullet,\bullet}(A)}
\newcommand{\CPAA}{C_{\Pe}^{\bullet,\bullet}(A)}
\newcommand{\emp}{\emptyset}
\newcommand{\trh}{\bar{C}^{\bullet}(A)}
\newcommand{\mi}{\left(\begin{array}{c}m\\i\end{array}\right)}

\newtheorem{lemma}{Lemma}

\newtheorem{rem}{Remark}
\newtheorem{ex}{Example}
\newtheorem{defn}{Definition}
\newtheorem{prop}{Proposition}

\begin{document}

\title{{\bf A master identity for homotopy Gerstenhaber algebras}}
\author{F\"{U}SUN AKMAN \\ Dept. of Mathematics and Statistics \\ Utah
State University \\ Logan, UT 84322-3900 \\ fusun@math.usu.edu}
\date{October 1, 1997}
\maketitle

\begin{abstract}
We produce a master identity $\tmtm$ for homotopy Gerstenhaber algebras, as
defined by Getzler and Jones and utilized by Kimura, Voronov, and Zuckerman
in the context of topological conformal field theories. To this end, we
introduce the notion of a ``partitioned multilinear map'' and explain the
mechanics of composing such maps. In addition, many new examples of pre-Lie
algebras and homotopy Gerstenhaber algebras are given.
\end{abstract}

\tableofcontents

\section{Introduction}

The proliferation of algebraic identities in string theory, such as the
lower identities for BRST complexes and topological vertex operator
algebras (TVOA) in Lian and Zuckerman's leading work \cite{LZ} (see also
\cite{PS}) and those given by Kimura, Voronov, and Zuckerman for homotopy
Gerstenhaber algebras -specifically for TVOA's- in \cite{KVZ}, has led the
author to study a common language for multilinear maps and their
compositions. Gerstenhaber's braces
\[ x\{ y\} =x\circ y \]
(which we will write as $\{ x\}\{ y\}$, following the notation of
\cite{KVZ}) have been in existence for more than thirty years and denote a
certain rule of composition of two multilinear maps $x$, $y$ on a vector
space $A$ with values in $A$. They were extended by Getzler in \cite{Get}
to the composition 
\[ x\{ y_1,\dots,y_n\} =\{ x\}\{ y_1,\dots,y_n\}\]
of several multilinear maps with one. These composition rules were
essential in redefining, and exploring new properties and examples of, 
commonly studied 
algebras like associative, Lie, homotopy associative ($\Ai$), homotopy Lie
($\Li$), Batalin-Vilkovisky, Gerstenhaber, and homotopy Gerstenhaber
($\Gi$) algebras. For example, many fundamental (sets of) identities, most
notably associativity, can be grouped under
\[ m\circ m=0,\]
where $m$ might be a homogeneous multilinear map or a formal infinite sum
of such maps. An overview will be given in Section~\ref{mubr}.
\vsp

In turn, the iterated composition operations have been shown to satisfy
many algebraic identities (see \cite{VG}) resembling those which arise in
the context of operads and homotopy algebras. Curiously, the multibraces
\[ \{ v_1,\dots,v_m\}\cdots\{ w_1,\dots,w_n\}\]
of \cite{KVZ}, defined on a TVOA or any $\Gi$ algebra, satisfy similar (but
more numerous) identities. As a first step in understanding the
similarities, the Gerstenhaber-Getzler braces were extended by the author
in \cite{A} to multiple substitutions/compositions of the form
\[ \{ x\}\{ y_1,\dots,y_m\}\cdots\{ z_1,\dots,z_n\}\]
on the Hochschild complex
\[ \CA =Hom(TA;A)\]
of multilinear maps on $A$, and eventually to
\[ \{ x_1,\dots,x_m\}\cdots\{ y_1,\dots,y_n\}\]
on
\[ \CAA =Hom(TA;TA)\]
(examples can be found in Section~\ref{mubr}). Although many new results
and identities were unveiled, the ``unadorned'' multilinear maps were
apparently not enough to describe the rich structure in \cite{KVZ} in a
nutshell.
\vsp

In this paper we introduce ``partitioned multilinear maps'' $x(\pi)$ which
are no different from ordinary multilinear maps except in composition. The
ordered partitions 
\[ \pi =\I\]
of nonnegative integers $i_1+\cdots +i_r$ give us a grouping of the
arguments of the map $x$, and the algebra generated by the ``products'' of
such partitions governs the composition rule: if 
\[ \pi\ast\pi'=\pi_1+\cdots +\pi_n,\]
then
\[ \{ x(\pi)\}\{ y(\pi')\} =\sum_{i=1}^nz_i(\pi_i),\]
where the resulting partitioned maps $z_i(\pi_i)$ are rigorously defined in
Sections \ref{ikub} and \ref{pmiz}. The product~$\ast$ is reduced to the 
familiar rule
\[ (i)\ast(j)=(i+j-1)\]
for singletons, which says that 
\[ \{ x(i)\}\{ y(j)\} =z(i+j-1),\]
or the composition of $i$-linear and $j$-linear maps is an
$(i+j-1)$-linear map. With this notation, it is possible to write the
algebraic master equation for homotopy Gerstenhaber algebras (which were 
previously defined only in terms of a topological operad) as
\[ \tmtm,\]
where
\[ m=\sum_{\pi}m(\pi)\]
is a formal sum of all the partitioned multilinear maps involved, and
$\tilde{m}$ is a term-by-term modification of~$m$ (only in
$\pm$~signs). This identity is to be interpreted as follows: the finite sum
of multilinear maps of ``type''~$\tip$, which result from the compositions
of all $m(\pi)$ and $m(\pi')$ such that
\[ \pi\ast\pi'=\tip +\cdots,\]
is identically zero for all ordered partitions $\tip$. If only singletons
$(i)$ are allowed as valid partitions, we obtain an $\Ai$ algebra. Some
interesting subalgebras of the $\Gi$ algebra are studied in
Section~\ref{subs}. A word of caution: the master identity may be
modified as
\[ \{\tilde{m}\}\{\tilde{m}\} +F(\tilde{m})=0;\]
see Remark~\ref{remar}.
\vsp

In addition to the master identity, which frees us from the necessity of
drawing several pictures for every given partition and supplies us with an
algorithm for writing the subidentities, we will introduce a variety of new
examples of homotopy Gerstenhaber algebras and pre-Lie algebras. For
example, the algebra $\Pe$ of regular partitions is a right pre-Lie algebra,
any vertex operator algebra with the Wick product is a left pre-Lie
algebra, and it may be possible to build up a $\Gi$ algebra from scratch 
starting with a square-zero Batalin-Vilkovisky type operator and the
$\Phi$-operators introduced by the author in \cite{ABV}. The
multibraces notation of \cite{KVZ} will be replaced by the partitioned map
notation $m\I$ so that braces can be used to denote composition. 
\vsp\vsp

{\it Acknowledgments.} I am very much indebted to Sasha Voronov, who
explained the pictures in \cite{KVZ} to me in great detail, and to Jim
Stasheff, for his continuing guidance and support. Many thanks are due
Haisheng Li and Chongying Dong for discussions on pre-Lie algebras and the
Wick and Zhu brackets.

\section{Partitioned maps}

\subsection{Multibraces}\label{mubr}

This is a short review of the multibraces notation introduced by
Gerstenhaber \cite{Ger} and further developed by Getzler \cite{Get}
and the author \cite{A}. 
For more details and historical references the reader is referred to \cite{A}.
\vsp

For the time being it suffices to consider compositions of
multilinear maps 
\[ x:A^{\ot n}\ra A\]
where $A$ is a (super) graded vector space
\[ A=\op_{n\in{\rm Z}}A^n.\]
These maps live in the {\bf Hochschild complex}
\[ \CA=\op_{n\geq 0}Hom(A^{\ot n};A),\]
which may be replaced by its ``completion'' 
\[ \CA =Hom(TA;A)\]
to accommodate formal infinite sums ($TA$ is the tensor algebra on
$A$). Elements of $A$, or maps from the (characteristic zero) field into
$A$, are treated on the same footing as higher multilinear maps. If $x$ is
an $n$-linear map ($n\geq 0$) which changes the total super degree of its
arguments by the integer amount $s$, it is assigned the gradings
\[ d(x)=n-1\;\;\;\;\;\;{\rm and}\;\;\;\;\;\; |x|=s.\]

The compact expression
\be \label{comp}\{ x\}\{ y_1,\dots,y_k\}\cdots\{ z_1,\dots,z_l\}\{
a_1,\dots,a_n\},\ee
where $x$, $y_i$, $z_j$ are maps and $a_1$,...,$a_n\in A$, indicates that:
\begin{itemize}
\item Every symbol except $x$ is to be substituted into one on the left, in
every possible way. But:
\item Symbols within the same pair of braces cannot be substituted into one
another, and the internal order within each pair must be retained. Then:
\item All such expressions are to be added up, each term accompanied by a
product of signs 
\[ (-1)^{d(u)d(v)+|u||v|}\]
whenever two symbols $u$, $v$ are interchanged with respect to (\ref{comp}).  
\end{itemize}
Super and $d$-degrees are preserved under such multiple
composition-substitutions. 

\begin{ex} \label{exam}
Let $d(x)=2$, $d(y)=1$, and $d(z)=2$. Then
\bea &&\{ x\}\{ y,z\}\{ a,b,c,d,e,f\}\nn\\
&=&(-1)^{d(z)(d(a)+d(b))+|z|(|a|+|b|)}x(y(a,b),z(c,d,e),f)\nn\\
&&+(-1)^{d(z)(d(a)+d(b)+d(c))+|z|(|a|+|b|+|c|)}x(y(a,b),c,z(d,e,f))\nn\\
&&+(-1)^{d(y)d(a)+|y||a|+d(z)(d(a)+d(b)+d(c))+|z|(|a|+|b|+|c|)}
x(a,y(b,c),z(d,e,f))\nn\\
&=&(-1)^{|z|(|a|+|b|)}x(y(a,b),z(c,d,e),f)\nn\\
&&+(-1)^{|z|(|a|+|b|+|c|)}x(y(a,b),c,z(d,e,f))\nn\\
&&-(-1)^{|y||a|+|z|(|a|+|b|+|c|)}x(a,y(b,c),z(d,e,f)).\nn\eea
But note that
\bea &&\{ x\}\{ y\}\{ z\}\{ a,b,c,d,e,f\}\nn\\
&=&\pm x(y(a,b),z(c,d,e),f)\pm x(y(a,b),c,z(d,e,f))\pm x(a,y(b,c),z(d,e,f))
\nn\\
&&\pm x(z(a,b,c),y(d,e),f)\pm x(z(a,b,c),d,y(e,f))\pm x(a,z(b,c,d),y(e,f))
\nn\\
&&\pm x(y(z(a,b,c),d),e,f)\pm x(y(a,z(b,c,d)),e,f)\pm x(a,y(z(b,c,d),e),f)
\nn\\
&&\pm x(a,y(b,z(c,d,e)),f)\pm x(a,b,y(z(c,d,e),f))\pm x(a,b,y(c,z(d,e,f))).
\nn\eea
\end{ex}

\begin{ex} By definition,
\[ \{ x\}\{ a_1\}\cdots\{ a_n\}=\sum_{\sigma\in S_n}(-1)^{{\rm sgn}(\sigma)}
x(a_{\sigma(1)},\dots,a_{\sigma(n)})\]
(super degrees ignored for simplicity).
\end{ex}

\begin{ex} If $m$ is an even bilinear map, the condition
\[ \{ m\}\{ m\} =0\]
is equivalent to associativity:
\bea && \{ m\}\{ m\}\{ a,b,c\}\nn\\
&=& m(m(a,b),c)+(-1)^{d(m)d(a)+|m||a|}m(a,m(b,c))\nn\\
&=& m(m(a,b),c)-m(a,m(b,c))=0.\nn\eea
\end{ex}

\begin{ex} For an even bilinear associative map $m$, the bracket defined by
\[ {[a,b]}_m =\{ m\}\{ a\}\{ b\} =m(a,b)-(-1)^{|a||b|}m(b,a)\]
is a (graded) Lie bracket, because
\[ \{ m\}\{ m\}\{ a,b,c\} =0\;\;\;\forall a,b,c\]
implies that
\[ \{ m\}\{ m\}\{ a\}\{ b\}\{ c\} =0\;\;\;\forall a,b,c,\]
which is equivalent to the Jacobi identity for ${[\; ,\;]}_m$.
\end{ex}

\begin{ex} A vector space $A$ with a bilinear multiplication map
\[ m(a,b)=a\star b\]
is a {\bf right pre-Lie algebra} if the identity
\be\label{prelie} (a\star b)\star c-a\star(b\star c)=
(-1)^{|b||c|}((a\star c)\star b-a\star(c\star b))\ee
holds (there may be any number of gradings on $A$). 
Gerstenhaber \cite{Ger} showed that his bigraded composition product
\[ x\circ y=\{ x\}\{ y\}\]
on -truncated- $\CA$ is a right pre-Lie product, 
and therefore the bigraded {\bf Gerstenhaber bracket} 
\[ {[x,y]}\defi\{ x\}\{ y\} -(-1)^{d(x)d(y)+|x||y|}\{ y\}\{ x\}\]
leads to a bigraded Lie algebra. This last statement is easy to see
in general for $(A,m)$ because (\ref{prelie}) can be expressed as
\[ \{ m\}\{ m\}\{ a,\{ b\}\{ c\}\,\} =0\;\;\;\forall a,b,c,\]
which implies
\[ \{ m\}\{ m\}\{ a\}\{ b\}\{ c\} =0\;\;\;\forall a,b,c.\]
Additional results about $\CA$ can be found in Sections \ref{ikud} and 
\ref{hrev}.
\end{ex}

\begin{ex}
In order to write master identities for {\bf homotopy associative and homotopy
Lie ($\Ai$ and $\Li$) algebras}, we are forced to introduce the single
grading
\[ ||x||=d(x)+|x| \]
on $\CA$, and to define
\[\tilde{x}(a_1,\dots,a_n)=(-1)^{(n-1)||a_1||+(n-2)||a_2||+\cdots
||a_{n-1}||}\{ x\}\{ a_1,\dots,a_n\}\]
for $d(x)=n-1$. Then, for example, an $\Ai$ algebra is nothing but a graded
vector space $A$ with a formal infinite sum of maps
\[ m=m(1)+m(2)+\cdots\]
satisfying
\[\{ \tilde{m}\}\{\tilde{m}\} =0.\]
In this setting we have 
\[ d(m(i))=i-1\;\;\;\;\;\;{\rm and}\;\;\;\;\;\; |m(i)|\equiv i\;\;\;
\mbox{(mod 2)},\]
so that $||m(i)||$ is always odd.
Moreover, an interchange of $u$ and $v$ results in the sign
\[ (-1)^{||u||\, ||v||}.\]
The master identity above reduces to the associativity condition on $m(2)$
when $m=m(2)$.
\end{ex}

\begin{ex} It is well-known \cite{LS} that graded antisymmetrization of the
products in an $\Ai$ algebra leads to an $\Li$ algebra. The proof in~\cite{A}
uses the following simple fact: the usual $\Li$ identities are equivalent
to
\[ \{ \tilde{m}\}\{ \tilde{m}\}\{ a_1\}\cdots\{ a_n\} =0\;\;\;\forall 
a_1,\dots,a_n\]
in this case, which is an immediate consequence of
\[ \{ \tilde{m}\}\{\tilde{m}\}\{ a_1,\dots,a_n\} =0\;\;\;\forall
a_1,\dots,a_n.\]
\end{ex}

\begin{ex} \label{sekiz}
In \cite{A}, the Hochschild complex was extended to
\[ \CAA =Hom(TA;TA)\]
by allowing multilinear maps with values in $TA$, including elements of
$TA$. The rules remain the same, but many formerly inadmissible expressions
have a meaning in the new complex. For example,
\[ \{ a_1,\dots,a_n\}\defi a_1\ot\cdots\ot a_n,\]
and
\[ \{ a\}\{ b\} =\{ a,b\}\pm\{ b,a\} =a\ot b-(-1)^{|a||b|}b\ot a.\]
\end{ex}

\subsection{The algebra of ordered partitions}

\subsubsection{Regular partitions}

We will define a binary multiplication $\ast$ on the free abelian group $\Pe$
spanned by the {\bf (regular) ordered partitions} 
\[ \pi=\I\;\;\;\;\;\; r\geq 1,\; i_1,\dots,i_r\geq 1\]
of all positive integers. These partitions  are to be thought of as
grouping the arguments of multilinear maps
\[ x:A^{\ot(i_1+\cdots +i_r)}\ra A\]
of {\bf type} $\pi$. We grade the basis elements $\pi$ of $\Pe$ by
\[ d(\pi)=i_1+\cdots +i_r-1\]
(total number of {\bf arguments} minus one) and
\[ \bar{d}(\pi)=r-1\]
(total number of {\bf slots} minus one), and remark that both degrees will
be preserved under multiplication. The product $\pi\ast\pi'$ of two ordered
partitions will be a sum
\[ \pi_1+\cdots +\pi_n \]
of basis elements, where the multiplicity of each $\pi_l$ is chosen to be 1
for convenience. We will later interpret the outcome as the collection of
types of multilinear maps that arise from the composition of two maps of
types $\pi$ and $\pi'$. An arbitrary element of the algebra $\Pe$ will be
denoted by $p$, $\hat{p}$, etc. {\bf Notation:}
If $\pi$ appears with a nonzero (integral) 
coefficient in $p$, we will write $\pi\ra p$.
\vsp

We first describe the product
\[ p=(i)\ast\J\]
with
\[ d(p)=i+j_1+\cdots +j_s-2,\;\;\;\bar{d}(p)=s-1\]
by
\be\label{bir} 
{p}=\sum_{u_1+\cdots +u_s=i-1,u_l\geq 0}(j_1+u_1|\cdots |j_s+u_s).\ee
The most general definition of 
\[ p=\pi\ast\pi'=\I\ast\J\]
is first approximated by
\be \hat{p}=\sum_{l=1}^r\sum_{\tip\ra(i_l)\ast\pi'}(i_1|\cdots |i_{l-1}|\tip |
i_{l+1}|\cdots i_r),\label{iki}\ee
after which we combine like terms under coefficient~1 and write
\be\label{uc} p=\sum_{\tip\ra\hat{p}}\tip,\ee
that is, we multiply one slot of $\pi$ with $\pi'$ at a time and insert the
resulting basis elements into the slot in question one by one, finally
erasing duplicate terms.

\begin{ex} In the following computation, we reduce the coefficient of
$(1|1|3|4)$ to~1 although it shows up twice in $\hat{p}$:
\[ (1|1|4)\ast(1|3)=(1|3|1|4)+(1|1|3|4)+(1|1|4|3)+(1|1|2|5)+(1|1|1|6).\]
\end{ex}

\begin{ex}
The special case
\[ (i)\ast(j)=(i+j-1)=(j)\ast(i)\]
reflects the fact that the composition of an $i$-linear map with a $j$-linear
map in the sense of Gerstenhaber is an $(i+j-1)$-linear map.
\end{ex}

The partition $(1)$ acts as a two-sided identity in the noncommutative
algebra $\Pe$. The product~$\ast$ is not associative nor left pre-Lie, but
it does turn out to be right pre-Lie.
These statements are false if we count multiplicities.

\begin{prop} The algebra $\Pe$ of ordered partitions is a right pre-Lie 
algebra under the product~$\ast$, namely, the identity
\[ (p_1\ast p_2)\ast p_3-p_1\ast(p_2\ast p_3)=(p_1\ast p_3)\ast p_2-
p_1\ast(p_3\ast p_2)\]
holds.
As a result, the (ungraded) Gerstenhaber bracket
\[ {[p_1,p_2]}\defi p_1\ast p_2-p_2\ast p_1 \]
on $\Pe$ is a Lie bracket.
\end{prop}

{\it Proof.} It suffices to give a proof for basis elements
\[ \pi_1=\I,\;\;\;\pi_2=\J,\;\;\;\pi_3=\K .\]

(A) Special case: we will show that when $\pi_1=(i)$ associativity holds on
the nose, i.e. 
\[ (\,(i)\ast\pi_2)\ast\pi_3=(i)\ast(\pi_2\ast\pi_3).\]
Let us compare the two sides (up to repetition of terms, hence the symbol 
$\approx$ instead of $=$) from the definition. We have
\[ (i)\ast\pi_2=\sum_{\usi}(j_1+u_1|\cdots |j_s+u_s)\]
and
\bea \label{Dag}&&(\,(i)\ast\pi_2)\ast\pi_3\\
&\approx&\sum_{\alpha =1}^s\sum_{\beta =1}^t\sum_{\usi}\sum_{\vtj}\nn\\
&&(j_1+u_1|\cdots |j_{\alpha -1}+u_{\alpha -1}|k_1+v_1|\cdots |k_t+v_t|
j_{\alpha +1}+u_{\alpha +1}|\cdots |j_s+u_s).\nn\eea
On the other hand, we have
\bea &&\pi_2\ast\pi_3\nn\\
&\approx&\sum_{\alpha =1}^s(j_1|\cdots |(j_{\alpha})\ast\pi_3|\cdots |
j_s)\nn\\
&\approx&\sum_{\alpha =1}^s\sum_{\beta =1}^t\sum_{\wtj}(j_1|\cdots |
j_{\alpha -1}|k_1+w_1|\cdots |k_t+w_t|j_{\alpha +1}|\cdots |j_s)\nn\eea
and
\bea\label{Dagg}&&(i)\ast(\pi_2\ast\pi_3)\\
&\approx&\sum_{\alpha =1}^s\sum_{\beta =1}^t\sum_{\wtj}\sum_{\zsi}\nn\\
&&(j_1+z_1|\cdots |j_{\alpha -1}+z_{\alpha -1}|k_1+w_1+z_{\alpha}|\cdots
|k_t+w_t+z_{\alpha +t-1}|j_{\alpha +1}+z_{\alpha +t}|\cdots
|j_s+z_{s+t-1}). \nn\eea
After making a (not necessarily invertible) change of variables
\bea &&u_1=z_1,\dots,u_{\alpha -1}=z_{\alpha -1}\nn\\
&&u_{\alpha}=z_{\alpha}+\cdots +z_{\alpha +t-1}\nn\\
&&u_{\alpha +1}=z_{\alpha +t},\dots,u_s=z_{s+t-1}\nn\\
&&v_1=w_1+z_{\alpha},\dots,v_t=w_t+z_{\alpha +t-1},\nn\eea
we obtain (\ref{Dag}) from (\ref{Dagg}). Note that any excess terms are
duplicates and can be discarded.
\vsp

(B) General case: we have
\[ \pi_1\ast\pi_2\approx\sum_{m=1}^r(i_1|\cdots |(i_m)\ast\pi_2|\cdots |
i_r)\]
and 
\bea &&(\pi_1\ast\pi_2)\ast\pi_3\nn\\
&\approx&\sum_{m=1}^r\sum_{\alpha =1}^{m-1}(i_1|\cdots |(i_{\alpha})\ast\pi_3|
\cdots |(i_m)\ast\pi_2|\cdots |i_r)\nn\\
&&+\sum_{m=1}^r(i_1|\cdots |(\,(i_m)\ast\pi_2)\ast\pi_3|\cdots |i_r)\nn\\
&&+\sum_{m=1}^r\sum_{\alpha =m+1}^r(i_1|\cdots |(i_m)\ast\pi_2|\cdots |
(i_{\alpha})\ast\pi_3|\cdots |i_r).\nn\eea
But then $\pi_1\ast(\pi_2\ast\pi_3)$ is just
\[ \pi_1\ast(\pi_2\ast\pi_3)=\sum_{m=1}^r(i_1|\cdots |(i_m)\ast
(\pi_2\ast\pi_3)|\cdots i_r),\]
and
\be\label{AA} (\pi_1\ast\pi_2)\ast\pi_3-\pi_1\ast(\pi_2\ast\pi_3)\ee
consists of all terms in which $\pi_2$ and $\pi_3$ (in any order!) are
multiplied on the right by two {\sl distinct} slots in $\pi_1$, thanks to
the special case. Obviously,
\be\label{BB} (\pi_1\ast\pi_3)\ast\pi_2-\pi_1\ast(\pi_3\ast\pi_2)\ee
consists of the exact same terms by symmetry, and the difference of 
(\ref{AA}) and (\ref{BB}) is zero.$\Box$

\begin{rem} Many right pre-Lie proofs follow the same pattern: one shows 
\[ (a\star b)\star c-a\star(b\star c)\]
is (possibly graded) symmetric in $b$ and $c$.
\end{rem}

\subsubsection{Partitions involving zeros}

In order to write down the $\Gi$ identities in a uniform fashion, we will
need to deal with partitions
\[ \pi=\I\]
with $i_1,\dots,i_r\geq 0$. Let us denote the basis of $\Pe$ consisting of
the regular partitions by $\Be$, the set of partitions with at least one
zero slot by $\Be_0$, their union $\Be\cup\Be_0$ by $\Beb$, and the free
abelian group spanned by $\Beb$ by $\Peb$. We extend the
multiplication~$\ast$ to $\Peb$ by the same rules
(\ref{bir})-(\ref{uc}). In effect, everything is as before except when there
is a zero slot on the {\sl left}: we have 
\[ (0)\ast(j_1|\cdots |j_s)=0 \]
since there are no nonzero numbers $u_1,\dots,u_s$ adding up to $-1$ (an
empty sum is zero).

\begin{ex}
\bea && (1|0|3)\ast(2)=(2|0|3)+(1|0|4)\nn\\
&& (2)\ast(1|0|3)=(1|0|4)+(1|1|3)+(2|0|3).\nn\eea
\end{ex}

\begin{rem}
The partition $(1)$ still acts as the left identity, but the algebra
$(\Peb,\ast)$ is not right pre-Lie any more. For example,
\[ (\,(i)\ast(0)\,)\ast(j)-(i)\ast(\,(0)\ast(j)\,)=(i+j-2),\]
whereas
\[ (\,(i)\ast(j)\,)\ast(0)-(i)\ast(\,(j)\ast(0)\,)=(i+j-2)-(i+j-2)=0\]
for $i\geq 2$, $j\geq 1$.
\end{rem}

\subsection{Partitioned multilinear maps}

\subsubsection{Regular partitioned maps} \label{ikub}

We introduce the notion of a {\bf (regular) partitioned multilinear map} 
\[ x(\pi)=x\I :A^{\ot i_1}\ot\cdots\ot A^{\ot i_r}\ra A,
\;\;\; \pi\in\Be \]
as an $(i_1+\cdots +i_r)$-linear map which is labelled by a regular
ordered partition
$\pi$ and is distinguished from ordinary maps only by the way we compose it
(Getzler and Jones mention multilinear maps $m_{k,l}$ which are similar to our
$m(k|l)$ in \cite{GJ2} but do not elaborate on their composition properties).
For substitution into $x(\pi)$, we will use the notation
\bea && \{ x\I\}\{ a_1,\dots,a_{i_1}|a_{i_1+1},\dots,a_{i_1+i_2}|\cdots |
a_{i_1+\cdots i_{r-1}+1},\dots,a_{i_1+\cdots +i_r}\}\nn\\
&=& \{ x\I\}\ar \nn\eea
rather than
\[ \{ x\I\}\{ a_1,\dots,a_{i_1+\cdots +i_r}\} .\]
When
\[ x=y_1(\pi_1)+\cdots +y_k(\pi_k)\;\;\;\;\;\;
\mbox{(with $d(\pi_1)=\cdots =d(\pi_k)=n-1$)} \]
is a sum of maps of different types, we will still write 
\[ \{ x\}\{ a_1,\dots,a_n\}\]
to mean
\[ \{ y_1(\pi_1)\}\ar +\cdots\;\;\;\;\;\;
\mbox{(say for $\pi_1=\I,\dots)$,}\]
where each map $y_l(\pi_l)$ is followed by the suitable barred braces. In
particular, we may use ordinary braces after the composition of two
maps $x(\pi)$ and $y(\pi')$ to denote the sum of all cases.
Recalling from the previous section that 
$\{ x(\pi)\}\{ y(\pi')\}$ is designed to be of the form
\[ \{ x(\pi)\}\{ y(\pi')\} =\sum_{\tip\ra\pi\ast\pi'}\zip,\]
we proceed to define each component $\zip$ rigorously. In the following 
account, the $\pm$~signs preceding each term are computed as in 
Section~\ref{mubr}, 
depending solely on the super degrees and the $d$-degrees of the maps and
elements involved; we will in general omit the full expressions to save
space. Note that
\[ d(x(\pi))=d(\pi)\;\;\;{\rm and}\;\;\;\bar{d}(x(\pi))=\bar{d}(\pi)\]
by definition.
\vsp

As with partitions, we start with the simpler case
\[ \{ x(i)\}\{ y(\pi')\},\;\;\;\pi'=\J.\]
Let
\[ \tip\ra(i)\ast(\pi'),\]
with 
\[ \tip=(j_1+u_1|\cdots |j_s+u_s)=(k_1|\cdots |k_s),\;\;\;u_1+\cdots +u_s=i-1\]
(note that all such $\tip$ are necessarily distinct). To define
\[ \{\zip\}\as,\]
we consider all possible subdivisions $S$ 
\bea \label{S} \{ a^{(1)}\} &=& \{ b^{(1)},c^{(1)},d^{(1)}\}\nn\\
&\vdots& \label{subdiv} \\
\{ a^{(s)}\} &=& \{ b^{(s)},c^{(s)},d^{(s)}\}\nn\eea
of the elements
\[ \{ a_1,\dots,a_{k_1+\cdots k_s}\} =\as \]
such that $c^{(1)}$ contains $j_1$ consecutive elements in $a^{(1)}$, 
$c^{(2)}$ has $j_2$ in $a^{(2)}$, and so on. Then we have
\bea && \{\zip\}\as \nn\\ &\defi& \sum_{S}\pm\{ x(i)\}\{ \, \{ b^{(1)}\}\cdots
\{ b^{(s)}\},\{ y\J\}\{ c^{(1)}|\cdots |c^{(s)}\},\{ d^{(1)}\}\cdots
\{ d^{(s)}\}\,\} .\nn\eea
The notation $\{ b^{(1)}\}\cdots\{ b^{(s)}\}$, for example, indicates that
we sum over all possible permutations of the remaining elements that come
{\sl before} the chosen strings $c^{(l)}$, while retaining the order 
within each $b^{(l)}$. The strings $b^{(l)}$ and $d^{(l)}$ may be empty.

\begin{ex} Since
\[ (3)\ast(2|4)=(2|6)+(3|5)+(4|4),\]
the composition $\{ x(3)\}\{ y(2|4)\}$ is a sum of three partitioned
maps. The first one is given by
\bea &&\{ x(3)\}\{ y(2|4)\}\{ a,b|c,d,e,f,g,h\}\nn\\
&=& \pm\{ x(3)\}\{\,\{y(2|4)\}\{ a,b|c,d,e,f\},g,h\}\nn\\
&&\pm\{ x(3)\}\{ c,\{ y(2|4)\}\{ a,b|d,e,f,g\},h\}\nn\\
&&\pm\{ x(3)\}\{ c,d,\{ y(2|4)\}\{ a,b|e,f,g,h\}\,\} ,\nn\eea
the second one by
\bea &&\{ x(3)\}\{ y(2|4)\}\{ a,b,c|d,e,f,g,h\}\nn\\
&=&\pm\{ x(3)\}\{\,\{ y(2|4)\}\{ a,b|d,e,f,g\},c,h\}\nn\\
&&\pm\{ x(3)\}\{\,\{ y(2|4)\}\{ a,b|d,e,f,g\},h,c\}\nn\\
&&\pm\{ x(3)\}\{ d,\{ y(2|4)\}\{ a,b|e,f,g,h\},c\}\nn\\
&&\pm\{ x(3)\}\{ a,\{ y(2|4)\}\{ b,c|d,e,f,g\},h\}\nn\\
&&\pm\{ x(3)\}\{ a,d,\{ y(2|4)\}\{ b,c|e,f,g,h\}\,\}\nn\\
&&\pm\{ x(3)\}\{ d,a,\{ y(2|4)\}\{ b,c|e,f,g,h\}\,\},\nn\eea
and the third one by
\bea &&\{ x(3)\}\{ y(2|4)\}\{ a,b,c,d|e,f,g,h\}\nn\\
&=&\pm\{ x(3)\}\{\,\{ y(2|4)\}\{ a,b|e,f,g,h\},c,d\}\nn\\
&&\pm\{ x(3)\}\{ a,\{ y(2|4)\}\{ b,c|e,f,g,h\},d\}\nn\\
&&\pm\{ x(3)\}\{ a,b,\{ y(2|4)\}\{ c,d|e,f,g,h\}\,\} .\nn\eea
\end{ex}

\begin{rem} Note that unlike the ordinary composition for non-partitioned
maps, the order of the elements $\as$ do change in the final
substitution. However, the order in each slot $a^{(l)}$ remains the same.
\end{rem}

In order to extend the definition to 
\[ \{ x(\pi)\}\{ y(\pi')\} =\{ x\I\}\{ y\J\},\]
we first identify the slot(s) in $\I$ which give rise to the partition
$\tip$ in the product $\pi\ast\pi'$ (when $r\geq 2$, it is possible to have
the outcome $\tip$ repeated in the multiplication process; we then reduce
the coefficient of $\tip$ in $\pi\ast\pi'$
to 1 and add all the multilinear maps of type~$\tip$ 
in $\{ x(\pi)\}\{ y(\pi')\}$). Without loss of generality, assume that the 
culprit is the first slot and that $\tip$ is not repeated. Then
\[ \tip=(j_1+k_1|\cdots |j_s+k_s|i_2|\cdots |i_r),\]
where
\[ k_1+\cdots +k_s=i_1-1,\]
and
\bea &&\{\zip\}\{ a^{(1)}|\cdots |a^{(s+r-1)}\}\nn\\
&\defi& \pm\{\,\{ x\I\}\{ -|a^{(s+1)}|\cdots |a^{(s+r-1)}\}\,\}\{ y\J\}\as .
\nn\eea
In other words, we fix the arguments of the remaining $r-1$ slots of $x(\pi)$,
and compose the two maps as if the first is an ordinary one with 
$\bar{d}=0$, according to the recipe in the previous paragraph.

\begin{ex} We verify that 
\[ (1|2)\ast(2|3)=(2|3|2)+(1|2|4)+(1|3|3),\]
and compute the first partitioned map which arises from multiplying the
first slot in $(1|2)$ by $(2|3)$:
\bea &&\{ x(1|2)\}\{ y(2|3)\}\{ a,b|c,d,e|f,g\}\nn\\
&=&\pm\{\,\{ x(1|2)\}\{ -|f,g\}\,\}\{ y(2|3)\}\{ a,b|c,d,e\} .\nn\eea
The last two maps come from multiplication with the second slot, hence
we have 
\bea &&\{ x(1|2)\}\{ y(2|3)\}\{ a|b,c|d,e,f,g\}\nn\\
&=&\pm\{\,\{ x(1|2)\}\{ a|-\}\,\}\{\,\{ y(2|3)\}\{ b,c|d,e,f\},g\}\nn\\
&&\pm\{\,\{ x(1|2)\}\{ a|-\}\,\}\{ d,\{ y(2|3)\}\{ b,c|e,f,g\}\,\},\nn\eea
and
\bea &&\{ x(1|2)\}\{ y(2|3)\}\{ a|b,c,d|e,f,g\}\nn\\
&=&\pm\{\,\{ x(1|2)\}\{ a|-\}\,\}\{\,\{ y(2|3)\}\{ b,c|e,f,g\},d\}\nn\\
&&\pm\{\,\{ x(1|2)\}\{ a|-\}\,\}\{ b,\{ y(2|3)\}\{ c,d|e,f,g\}\,\} .\nn\eea
\end{ex}

\begin{rem}[Shortcut] Suppose we are trying to compose $x(\pi)$ and
$y(\pi')$, with
\[ \pi=\I\;\;\;{\rm and}\;\;\;\pi'=\J .\]
Once we compute the product $\pi\ast\pi'$ and decide which 
\[ \tip =\K\ra\pi\ast\pi'\]
to go after, we look at $\att$ and check out which consecutive $s$~slots are
large enough to accommodate $y(\pi')$. The remaining elements in these
slots as well as $y$ will go into one slot of $x(\pi)$, therefore we also
have to make sure that the remaining slots and the combination of chosen
slots will conform to the type of $x$ before deciding this is one of the
feasible combinations.
\end{rem}

\begin{ex} In \cite{ABV}, we defined {\bf higher order differential operators}
$\De :A\ra A$ for a noncommutative, nonassociative algebra $(A,m)$ to
coincide with the commutative and associative case given by Koszul
in~\cite{Ko}. A follow-up on the simplification of notation due to multibraces,
and a  generalization, appeared in~\cite{A}. A linear operator $\De$ is
called a differential operator of order~$\leq r$ if and only if a certain
$(r+1)$-linear map
\be \Phi_{\De}^{r+1}(a_1,\dots,a_{r+1})\label{phiop}\ee
is identically zero; we now want to think of (\ref{phiop}) as a partitioned
map $\Phi_{\De}^{r+1}(r|1)$, because the last slot is distinguished and the
original inductive definition 
\bea \Phi_{\De}^1(a)&=&\De(a)\nn\\
\Phi_{\De}^2(a,b)&=&\De(ab)-\De(a)b-(-1)^{|\De ||a|}a\De(b)\nn\\
\Phi_{\De}^{r+2}(a_1,\dots,a_r,b,c)&=&\Phi_{\De}^{r+1}(a_1,\dots,a_r,bc) 
-\Phi_{\De}^{r+1}(a_1,\dots,a_r,b)c\nn\\
&&-(-1)^{|b|(|\De |+|a_1|+\cdots +|a_r|)}
b\Phi_{\De}^{r+1}(a_1,\dots,a_r,c)\;\;\; r\geq 1\nn\eea
in \cite{ABV} can be conveniently rewritten as
\bea \phi_{\De}^1(a)&=&\{\De\}\{ a\}\nn\\
\Phi_{\De}^2(a,b)&=&[\Phi_{\De}^1(1),m(2)]\{ a,b\}\nn\\
\Phi_{\De}^{r+2}(a_1,\dots,a_r,b,c)&=&[\Phi_{\De}^{r+1}(r|1),m(2)]
\{ a_1,\dots,a_r|b,c\}\;\;\; r\geq 1\nn\eea
in terms of Gerstenhaber brackets and partitioned maps. Note that
\[ (r|1)\ast(2)=(2)\ast(r|1)=(r+1|1)+(r|2).\]
Here we may designate the same multilinear map $\Phi_{\De}^{r+2}$ to be of
type $(r+2)$, $(r+1|1)$, or $(r|2)$, depending on what we want to
describe! More in Section~\ref{phiope}.
\end{ex}

\subsubsection{Higher products of regular partitions}\label{hprp}

The algebra $\Pe$ enjoys higher products
\be N(1|\la_1|\cdots |\la_t):\Pe\ot\Pe^{\ot\la_1}\ot\cdots\ot\Pe^{\ot\la_t}
\ra\Pe, \label{N}\ee
where $(\la_1|\cdots |\la_t)$ itself is any regular partition, and
\[ N(1|1)(\pi,\pi')\defi\pi\ast\pi'.\]
Such products will be used
to model higher compositions of partitioned maps in the next section.
\vsp

We had better define 
\[ \{ N(1|k)\}\{ (i)|(j_1),\dots,(j_k)\}\defi(i+j_1+\cdots +j_k-k),\]
because we will obtain an $(i+j_1+\cdots +j_k-k)$-linear map when we
compute $\{ x(i)\}\{ y_1(j_1),\dots,y_k(j_k)\}$: we simply add up all
$d$-gradings. Similarly, it makes sense to define
\[ \{ N(1|1|\cdots |1)\}\{ (i)|(j_1)|\cdots |(j_k)\}\defi(i+j_1+\cdots
+j_k-k),\]
because again
\[ d(\{ x\}\{ y_1\}\cdots\{ y_k\})=d(x)+d(y_1)+\cdots +d(y_k).\]
Note that $\bar{d}$ is preserved as well in both cases.
\vsp

Next, we define for $i\geq 2$
\bea &&\{ N(1|2)\}\{ (i)|\J,\K\}\nn\\
&\defi&\sum_{u_1+\cdots +u_{s+t-1}=i-2,u_l\geq 0}(j_1+u_1|\cdots |j_{s-1}
+u_{s-1}|j_s+k_1+u_s|k_2+u_{s+1}|\cdots |k_t+u_{s+t-1}),\nn\eea
that is, we write down the partition 
\[ (j_1|\cdots |j_s+k_1|\cdots |k_t),\]
and add to the slots nonnegative integers with total $i-2$ in every
possible way. In case of
\[ \{ N(1|3)\}\{ (i)|\J,\K,(l_1|\cdots |l_u)\}\;\;\;(i\geq 3),\]
we start with
\[ (j_1|\cdots |j_{s-1}|j_s+k_1|k_2|\cdots |k_{t-1}|k_t+l_1|l_2|\cdots 
|l_u),\]
and this time distribute a sum of $i-3$. The generalization to
\[ \{ N(1|k)\}\{ (i)|\pi_1,\dots,\pi_k\},\;\;\;\pi_l\in\Be\]
is clear. Once more, both $d$ and $\bar{d}$ are preserved. In the next
section we will see the definition of the composition $\{ x\}\{ y,z\}$
which motivates the above formulas.
\vsp

As for
\[ \{ N(1|1|\cdots |1)\}\{ (i)|\pi_1|\cdots |\pi_k\},\]
we recall the meaning of $\{ x\}\{ y\}\{ z\}$ for non-partitioned maps from
Section~\ref{mubr} (see Example~\ref{exam}):
\[ \{ x\}\{ y\}\{ z\} =\{ x\}\{ y,z\}\pm\{ x\}\{ z,y\}\pm\{ x\}\{\,\{ y\}
\{ z\}\,\} .\]
In other words, $y$ and $z$ can be substituted into $x$ separately (and in
any order), or $z$ can go into $y$ first. Then by analogy we want to have
\bea &&\{ N(1|1|1)\}\{ (i)|\pi_1|\pi_2\}\nn\\
&\approx&\{ N(1|2)\}\{ (i)|\pi_1,\pi_2\} +\{ N(1|2)\}\{ (i)|\pi_2,\pi_1
\} +\{ N(1|1))\}\{ (i)|\{ N(1|1)\}\{\pi_1|\pi_2\}\,\},\nn\eea 
where the symbol $\approx$ again means that any repetition of basis
elements $\tip$ on the right is to be ignored. The formula should hold even
when the first partition $(i)$ is replaced by any $\pi\in\Be$. In fact, all
higher products~(\ref{N}) should be constructed from lower products by
thinking in terms of the corresponding compositions.

\begin{ex} We have
\bea &&\{ x\}\{ y_1,y_2\}\{ z\}\nn\\
&=&\{ x\}\{ y_1,y_2,z\}\pm\{ x\}\{ y_1,z,y_2\}\pm\{ x\}\{ z,y_1,y_2\}\nn\\
&&\pm\{ x\}\{\,\{ y_1\}\{ z\},y_2\}\pm\{ x\}\{ y_1,\{ y_2\}\{ z\}\,\},\nn\eea
therefore by definition
\bea &&\{ N(1|2|1)\}\{\pi |\pi_1,\pi_2|\pi'\}\nn\\
&\approx&\{ N(1|3)\}\{\pi |\pi_1,\pi_2,\pi'\} +\{ N(1|3)\}\{\pi|\pi_1,\pi',
\pi_2\} +\{ N(1|3)\}\{\pi |\pi',\pi_1,\pi_2\}\nn\\
&&+\{ N(1|2)\}\{\pi |\{ N(1|1)\}\{\pi_1|\pi'\},\pi_2\} +\{ N(1|2)\}
\{\pi |\pi_1,\{ N(1|1)\}\{\pi_2 |\pi'\}\,\} .\nn\eea
\end{ex}

To complete the discussion of higher products, we replace $(i)$ by $\pi
=\I$: then by definition,
\[ N(1|\la_1|\cdots |\la_t)\}\{\pi |\cdots\}\approx\sum_{\alpha =1}^r
(i_1|\cdots |i_{\alpha -1}|\{ N(1|\la_1|\cdots
|\la_t)\}\{(i_{\alpha})|\cdots\} |i_{\alpha +1}|\cdots |i_r).\]

\begin{ex} 
\bea &&\{ N(1|2)\}\{(2|3)|(2),(3|4)\}\nn\\
&\approx&(\{ N(1|2)\}\{(2)|(2),(3|4)\} |3)+(2|\{ N(1|2)\}\{(3)|(2),
(3|4)\})\nn\\
&=&(5|4|3)+(2|5|5)+(2|6|4).\nn\eea
\end{ex}

\subsubsection{Higher compositions of regular partitioned maps}

Following \cite{Get} and \cite{A}, we would like to define 
\[ \{ x(\pi)\}\{ y_1(\pi_1),\dots,y_k(\pi_k)\},\]
and ultimately
\[ \{ x(\pi)\}\{ y_1(\pi_1),\dots,y_k(\pi_k)\}\cdots\{ z_1(\pi_1'),\dots,
z_l(\pi_l')\}\]
for $\pi$, $\pi_i$, $\pi_j'\in\Be$. The rules will be similar to the
nonpartitioned case, and types of maps will be governed by the enriched 
algebra of
partitions discussed in the previous section. To give an easy example, we
consider 
\[ \{ x(i)\}\{ y(\pi_1),z(\pi_2)\},\;\;\;{\rm with}\;\;\;\pi_1=\J,
\pi_2=\K\in\Be .\]
We choose
\[ \tip\ra\{ N(1|2)\}\{(i)|\pi_1,\pi_2\}\]
and set up variables
\[ \{ a^{(1)}|\cdots |a^{(s+t-1)}\}\]
accordingly. Then
\be\label{xyz}\{ x(i)\}\{ y(\pi_1),z(\pi_2)\}\{ a^{(1)}|\cdots |
a^{(s+t-1)}\}\ee
will be defined as a sum over subdivisions $S$ of these variables. Each $S$
will look like
\bea a^{(1)}&=&\{ b^{(1)},c^{(1)},d^{(1)}\}\nn\\ 
{} &\vdots& {}\nn\\
a^{(s-1)}&=&\{ b^{(s-1)},c^{(s-1)},d^{(s-1)}\}\nn\\
a^{(s)}&=&\{ b^{(s)},c^{(s)},d^{(s)}=e^{(1)},f^{(1)},g^{(1)}\}\nn\\
a^{(s+1)}&=&\{ e^{(2)},f^{(2)},g^{(2)}\}\nn\\
{} &\vdots& {}\nn\\
a^{(s+t-1)}&=&\{ e^{(t)},f^{(t)},g^{(t)}\},\nn\eea
where $c^{(l)}$ has $j_l$ consecutive elements, $f^{(l)}$ has $k_l$
consecutive elements, and the strings $b^{(l)}$, $d^{(l)}$, $e^{(l)}$, 
$g^{(l)}$ may be empty. Then the contribution of $S$ to~(\ref{xyz}) will be
\bea \pm\{ x(i)\}&&\{\,\{ b^{(1)}\}\cdots\{ b^{(s)}\},\{ y(\pi_1)\}\{
c^{(1)}|\cdots |c^{(s)}\},\{ d^{(1)}\}\cdots\{ d^{(s-1)}\}\{ e^{(1)}\}\cdots
\{ e^{(t)}\},\nn\\
&&\{ z(\pi_2)\}\{ f^{(1)}|\cdots |f^{(t)}\},\{ g^{(1)}\}\cdots\{ g^{(t)}\}
\,\} .\nn\eea

Any multiple compositions involving more than two pairs of braces can be
handled as sums and iterations of the latter: for example,
\[ \{ x(\pi)\}\{ y(\pi_1)\}\{ z(\pi_2)\}\]
means
\[ \{ x(\pi)\}\{ y(\pi_1),z(\pi_2)\}\pm\{ x(\pi)\}\{ z(\pi_2),y(\pi_1)\}\pm
\{ x(\pi)\}\{\,\{ y(\pi_1)\}\{ z(\pi_2)\}\,\} .\]

\subsubsection{Partitioned Hochschild complex}\label{ikud}

Let us denote the vector space of $n$-linear maps on $A$ ($n\geq 0$)
with values in $A$ by 
\[ Hom(A^{\ot n};A)\]
and the vector space of multilinear maps on $A$ of type $\pi=\I\in\Be
\cup\{ (0)\}$ by
\[ Hom(A^{\pi};A)=Hom(A^{\ot i_1}\ot\cdots\ot A^{\ot i_r};A)\]
(the two spaces are clearly isomorphic if $d(\pi)=n-1$). Depending on
whether we want infinite sums or not, we defined the regular Hochschild
complex $\CA$ to be either of the spaces 
\[ \op_{n\geq 0}Hom(A^{\ot n};A)\subset Hom(\op_{n\geq 0}A^{\ot n};A)=
Hom(TA;A).\]
In a similar manner, we may define the {\bf partitioned Hochschild complex}
$\CPA$ by either of
\[ \op_{\pi\in\Be\cup\{ (0)\}}Hom(A^{\pi};A)\subset Hom(\op_{\pi\in\Be
\cup\{ (0)\}}A^{\pi};A)=Hom(T(TA);A).\]

\begin{prop}[Gerstenhaber] The {\bf truncated Hochschild complex}
\[ \bar{C}^{\bullet}(A)=\op_{n\geq 1}Hom(A^{\ot n};A)\]
is a right pre-Lie algebra under the multiplication
\[ x\circ y=\{ x\}\{ y\} .\]\label{propo2}
\end{prop}

{\it Proof.} See \cite{Ger}, \cite{A}, and Proposition~\ref{dort} below.$\Box$
\vsp

The following Proposition answers a natural question:

\begin{prop}\label{dort} The {\bf truncated partitioned Hochschild complex}
\[ \bar{C}_{\Pe}^{\bullet}(A)=\op_{\pi\in\Be}Hom(A^{\pi};A)\]
is a right pre-Lie algebra under the composition of partitioned maps.
\end{prop}

{\it Proof.} The triple composition
\[ \{\,\{ x(\pi)\}\{ y(\pi')\}\,\}\{ z(\pi'')\} -\{ x(\pi)\}\{\,\{ y(\pi')
\}\{ z(\pi'')\}\,\}\]
consists of partitioned maps $w(\tip)$ obtained by
substituting $y$ and $z$ {\sl separately} into $x$ (we subtract the
terms where $z$ goes into $y$). Then by symmetry 
\[ \{\,\{ x(\pi)\}\{ z(\pi'')\}\,\}\{ y(\pi')\} -\{ x(\pi)\}\{\,\{ z(\pi'')
\}\{ y(\pi')\}\,\}\]
has the same summands up to an overall sign.$\Box$
\vsp

\begin{rem} The result holds for the full complex in both cases, keeping
in mind that $\{ a\}\{ b_1,\dots,b_n\} =0$ for $a$, $b_i\in A$.
\end{rem}

\begin{prop} The full {\bf extended Hochschild complex} 
\[ \CAA =Hom(\op_{i\geq 0}A^{\ot i};\op_{j\geq 0}A^{\ot j})\subset Hom(TA;TA)\]
defined in \cite{A} is a right pre-Lie algebra under the multiplication
\[ \{ x\}\{ y\} .\]
\end{prop}

{\it Proof.} Same as above. The difference is that now {\sl all} compositions 
$\{ x\}\{ y\}$ have a natural definition in the new complex, even when $x$
and $y$ are elements of $TA$ (see Example~\ref{sekiz} in
Section~\ref{mubr}).$\Box$

\subsubsection{Partitioned maps involving zeros}\label{pmiz}

We now need a theory of compositions of maps of type~$\pi$ for $\pi\in\Beb$,
consistent with our earlier definitions. It is clear that $x(0)$ (with
$d(x)=-1$) has to be an element of~$A$.

\begin{ex} For $n\geq 2$, we have
\[ (n)\ast(0)=(n-1),\]
and not surprisingly,
\[ \{ x(n)\}\{ y(0)\}\{ a_1,\dots,a_{n-1}\} =\{ x(n)\}\{ b\}\{
a_1,\dots,a_{n-1}\} \]
if $y(0)=b$.
\end{ex}

On the other hand,

\begin{ex} The expression $\{ x(0)\}\{ y(0)\}$ must be interpreted as zero
in $\CPA$, as in $\CA$ (after all, $(0)\ast(0)=0$). It can be thought of as
an element
\[ a\ot b-(-1)^{|a||b|}b\ot a\in A^{\ot 2}\]
inside
\[ \CAA =Hom(TA;TA)\]
(or its suitably defined counterpart $\CPAA$) only if we enlarge the
Hochschild complex. Similarly, the composition 
$\{ y(0)\}\{ x(n)\}$ for $n\geq 1$ is zero in $\CPA$.
\end{ex}

We need only define the subdivisions $S$ corresponding to some
\[ \{ x(i)\}\{ y\J\}\as \;\;\;\mbox{with $i\geq 1$.}\]
We again subdivide elements in each slot according to~(\ref{subdiv}), but
this time the number of elements in $c^{(l)}$ will be zero if
$j_l=0$. Since $b^{(l)}$ and $d^{(l)}$ were already allowed to be empty,
this is not a novel concept! 

\begin{ex} We have
\[ (4)\ast(1|0)=(1|3)+(2|2)+(3|1)+(4|0).\]
Then the six subdivisions for $(2|2)$, or $\{ a,b|c,d\}$, will be
\bea \emp,\, \{ a\},\, \{ b\} &|& \emp,\,\emp,\,\{ c,d\}\nn\\
\emp,\,\{ a\},\,\{ b\} &|&\{ c\},\,\emp,\,\{ d\}\nn\\
\emp,\,\{ a\},\{ b\} &|&\{ c,d\},\,\emp,\,\emp\nn\\
\{ a\},\,\{ b\},\,\emp &|&\emp,\,\emp,\,\{ c,d\}\nn\\
\{ a\},\,\{ b\},\,\emp &|&\{ c\},\,\emp,\,\{ d\}\nn\\
\{ a\},\,\{ b\},\,\emp &|&\{ c,d\},\,\emp,\,\emp,\nn\eea
and we will have
\bea &&\{ x(4)\}\{ y(1|0)\}\{ a,b|c,d\}\nn\\
&=&\{ x(4)\}\{\,\{ y(1|0)\}\{ a\},\{ b\}\{ c,d\}\,\}\nn\\
&&\pm\{ x(4)\}\{ c,\{ y(1|0)\}\{ a\},\{ b\}\{ d\}\,\}\nn\\
&&\pm\{ x(4)\}\{ c,d,\{ y(1|0)\}\{ a\},b\}\nn\\
&&\pm\{ x(4)\}\{ a,\{ y(1|0)\}\{ b\},c,d\}\nn\\
&&\pm\{ x(4)\}\{\,\{  a\}\{ c\},\{ y(1|0)\}\{ b\},d\}\nn\\
&&\pm\{ x(4)\}\{\,\{ a\}\{ c,d\},\{ y(1|0)\}\{ b\}\,\} \nn\eea
as the $(2|2)$-component of the composition.
\end{ex}

\begin{ex} For $i\geq 2$, the sum
\[ \sum_{j_1+\cdots +j_s=1}\{ x(i)\}\{ y\J\}\as\]
with each $y\J$ as the (linear) identity map on $A$ gives us exactly
\[ i\{ x(i)\}\{ a^{(1)}\}\cdots\{ a^{(s)}\} .\]
In particular,
\[\{ x(i)\}\{{\rm id}(1)\}\{ a_1,\dots,a_i\}=i\{ x(i)\}\{ a_1,\dots,a_i\} .\]
\label{duh} \end{ex}

\section{The master identity for homotopy Gerstenhaber algebras}

The present work on the master identity for $\Gi$ algebras originated from
a joint paper of Kimura, Voronov, and Zuckerman \cite{KVZ}, where a nontrivial
algebra over a certain cellular operad $\KM$ (first defined by Getzler and
Jones \cite{GJ}, based on ideas of Fox-Neuwirth) is described. Since 
\cite{KVZ} is an
excellent introduction to the subject and contains an in-depth description
of $\KM$, we will skip the origins and details of the following definition:

\begin{defn} A {\bf homotopy Gerstenhaber algebra ($\Gi$ algebra)} is an
algebra over the operad $\KM$.
\end{defn}

Instead, in this section, we will expand the discussion in \cite{KVZ} of 
how the individual $\Gi$ identities are obtained from certain
``configurations''. Our explanations will not require any background, and
the pictures will eventually coagulate into the equivalent

\begin{defn} \label{defiki}
A {\bf homotopy Gerstenhaber algebra ($\Gi$ algebra)} is a
super graded vector space
\[ A=\op_{n\in{\rm Z}}A^n\]
together with a collection
\[ m(\pi),\;\;\;\;\;\; |m(\pi)|\equiv d(\pi)+\bar{d}(\pi)+1\;\;\;
\mbox{(mod 2)}\]
of partitioned multilinear maps, 
such that 
\[ m(1|0)=m(0|1)={\rm id}_A \]
are the only nonzero ones among $m(\pi)$ with $\pi\in\Be_0$,
where the formal sum
\[ m=\sum_{\pi\in\Beb}m(\pi)\]
satisfies
\be\label{master}\{\tilde{m}\}\{\tilde{m}\} =0.\ee
\end{defn}

\begin{rem}\label{rema} Identity (\ref{master}) means the finite sum of
partitioned maps of type~$\tip$ is equal to zero in the composition for
each partition~$\tip$. It is understood that if a map~$m(\tip)$ is
identically zero, then we do NOT consider the $\tip$ part of~(\ref{master})
among the $\Gi$ subidentities. The Hochschild complex in Section~\ref{hrev}
is a good example. Also, we define the {\bf total degree} $||m(\pi)||$ to
be 
\[ ||m(\pi)||\defi d(\pi)+\bar{d}(\pi)+|m|,\]
which is always odd in a $\Gi$ algebra. The actual super degree comes from
dimensions of chains in the $\Gi$ operad $\KM$.
\end{rem}

We emphasize again that the elegant notation of \cite{KVZ} has to be
changed in
order to reserve the older notation of braces for compositions of
maps. The symbol $m\I$ will be used to denote the multilinear map in 
\cite{KVZ} shown by $r$ pairs of braces with $i_1,\dots,i_r$ arguments
respectively.
\vsp

The configurations which ``contribute to the differential'' and hence show
up in the $\Gi$ algebra identities are as follows: for each regular
ordered partition
\[ \tip=\K\in\Be \]
we draw $t$ vertical lines with $k_1,\dots,k_t$ points 
on them respectively from left to right, labelled by
\[ \{ a_1,\dots,a_{k_1+\cdots +k_t}\} =\att \]
(the lexicographical ordering is from
top to bottom and left to right). We then make a list of all possible
configurations in which one smooth ``bubble'' surrounds several points in
this picture, subject to the following conditions:
\begin{itemize}
\item If two points on the same vertical line are in the bubble, so are all
points between them.
\item If two points on different vertical lines are in the bubble, so is at
least one point from each line in between the two lines.
\end{itemize}
Note that bubbles around only one point and the giant bubble enclosing all
points are acceptable. We will call these configurations Type~I. Next, we
introduce bubbles that contain exactly one point on one line and an empty
portion of an adjacent line and call the resulting pictures Type~II. Bubbles
with empty portions positioned between different points will be counted 
as different (the empty portion may also be above or below all points on the
second line).
If we start with only one vertical line, there will be no Type~II pictures.
\vsp

Type~I pictures translate into the quadratic parts of the lower identities 
given
in \cite{KVZ} and those of Type~II into the seemingly non-quadratic parts. The
goal in the Type~I situation is to collect together all pictures which
depict the composition of $m(\pi)$ with $m(\pi')$ ($\pi,\pi'\in\Be$) 
such that $\tip\ra\pi\ast\pi'$. Each picture is really one subdivision $S$
for some $\{ m(\pi)\}\{ m(\pi')\}$. 
The original picture without the bubble represents 
\[ \{ m\K\}\att ;\]
we understand that each vertical line corresponds to one slot. Similarly,
the inside of a bubble is what we would have denoted in Section~\ref{ikub} by
\[ \{ m\J\}\{ c^{(1)}|\cdots |c^{(k)}\} \]
for a particular subdivision $S$ of $\att$, sitting as an element of~$A$
inside some $m\I$. Any
vertical lines seen above the bubble will have the leftover elements
which had formerly been denoted by 
\[ b^{(1)},\dots,b^{(k)}\]
and the ones below will have the elements
\[ d^{(1)},\dots,d^{(k)}.\]
We transcribe each of these pictures into symbols and add up! 
\vsp

In order to complete the $\tip$ portion of the identity $\{\tilde{m}\}\{ 
\tilde{m}\} =0$ we introduce a picture containing two vertical
lines with one point on the left (resp. right) and no points on the right 
(resp. left) to mean $m(1|0)$ (resp. $m(0|1)$). By Example~\ref{duh}
we know that the overall contribution of Type~II pictures for $\tip=\K$
will be
\bea \sum_{\alpha =1}^{t-1}(k_{\alpha}+k_{\alpha +1})&&\{ m(k_1|\cdots |
k_{\alpha -1}|k_{\alpha}+k_{\alpha +1}|k_{\alpha +2}|\cdots |k_t)\}\nn\\
&&\{ a^{(1)}|\cdots |a^{(\alpha -1)}|\{ a^{(\alpha)}\}\{ a^{(\alpha +1)}\} |
a^{(\alpha +2)}|\cdots |a^{(t)}\} .\nn\eea

Finally, we set the sum of all Type~I and Type~II terms equal to zero. 

\begin{ex}
For example, two vertical lines with points $\{ a\}$ and $\{ b,c\}$ 
respectively 
(i.e. the original picture for $\tip =(1|2)$) give us seven Type~II pictures,
which will add up to
\[ 3\{ m(3)\}\{ a\}\{ b,c\} .\]
Then the overall subidentity for $(1|2)$ is
\bea &&\{ m(1)\}\{ m(1|2)\}\{ a|b,c\}\pm\{ m(1|2)\}\{ m(1)\}\{ a|b,c\}\nn\\
&&\pm\{ m(2)\}\{\,\{ m(1|1)\}\{ a|b\},c\}\pm\{ m(2)\}\{ b,\{ m(1|1)\}\{
a|c\}\,\}\nn\\
&&\pm\{ m(1|1)\}\{\ a|\{ m(2)\}\{ b,c\}\,\}\pm 3\{ m(3)\}\{ a\}\{ b,c\} =0.
\nn\eea
This is Eqn.~(8) of \cite{KVZ} with the exception of the coefficient~3; see
Remark~\ref{remar} below. To check that all degrees are preserved, one must
take into account $\bar{d}(\{ a|b,c\})=1$ and $\bar{d}(\{ a\}\{ b,c\})=0$.
\end{ex}

\begin{ex} Let us write down all terms of $\tmtm$ which correspond to the 
partition $\tip =(1|2|3)$. We start with a picture having three vertical
lines, with points $\{ a\}$, $\{ b,c\}$, and 
$\{ d,e,f\}$ respectively, from left to right. Here is a list of all 
bubbles in Type~I pictures and the terms they contribute (up to sign):
\[ \begin{array}{|c|c|l|} 
\hline
\mbox{{\bf bubbles}} & \mbox{{\bf relevant product}} & \mbox{{\bf contributed 
terms}} \\
\mbox{{\bf in Type I pictures}} & \pi\ast\pi'=(1|2|3)+\cdots & \mbox{(up
to sign)} \\
\hline\hline
abcdef & (1)\ast(1|2|3) & \{ m(1)\}\{ m(1|2|3)\}\abc \\ \hline
a,\, b,\, c,\, d,\, e,\, f & (1|2|3)\ast(1) & \{ m(1|2|3)\}\{ m(1)\}\abc \\
\hline
abdef,\, acdef & (2)\ast(1|1|3) & \{ m(2)\}\{ m(1|1|3)\}\abc \\ \hline
bc & (1|1|3)\ast(2) & \{ m(1|1|3)\}\{ m(2)\}\abc \\ \hline
abcde,\, abcef & (2)\ast(1|2|2) & \{ m(2)\}\{ m(1|2|2)\}\abc \\ \hline
de,\, ef & (1|2|2)\ast(2) & \{ m(1|2|2)\}\{ m(2)\}\abc \\ \hline
abde,\, abef,\, acde,\, acef & (3)\ast(1|1|2) & \{ m(3)\}\{
m(1|1|2)\}\abc \\ \hline
abcd,\, abce,\, abcf & (3)\ast(1|2|1) & \{ m(3)\}\{ m(1|2|1)\}\abc \\ \hline
def & (1|2|1)\ast(3) & \{ m(1|2|1)\}\{ m(3)\}\abc \\ \hline
abd,\, abe,\, abf,\, acd,\, ace,\, acf & (4)\ast(1|1|1) & \{ m(4)\}
\{ m(1|1|1)\}\abc \\ \hline
bcdef & (1|1)\ast(2|3) & \{ m(1|1)\}\{ m(2|3)\}\abc \\ \hline
ab,\, ac & (2|3)\ast(1|1) & \{ m(2|3)\}\{ m(1|1)\}\abc \\ \hline
bdef,\, cdef & (1|2)\ast(1|3) & \{ m(1|2)\}\{ m(1|3)\}\abc \\ \hline
abc,\, bde,\, bef,\, cde,\, cef & (1|3)\ast(1|2) & \{ m(1|3)\}\{ m(1|2)\}
\abc \\ \hline
bcde,\, bcef & (1|2)\ast(2|2) & \{ m(1|2)\}\{ m(2|2)\}\abc \\ \hline
bcd,\, bce,\, bcf & (1|3)\ast(2|1) & \{ m(1|3)\}\{ m(2|1)\}\abc \\ \hline
bd,\, be,\, bf,\, cd,\, ce,\, cf & (1|4)\ast(1|1) & \{ m(1|4)\}\{ m(1|1)\}
\abc \\ \hline
\end{array}\]
It can be checked 
that all products $\pi\ast\pi'$ such that $(1|2|3)\ra\pi\ast\pi'$
are in the above list, and for each such product every subdivision $S$ 
leading to $(1|2|3)$
is given by a bubble! Next, we list the two Type~II pictures and their
contributions:
\[ \begin{array}{|c|c|l|}
\hline
\mbox{{\bf bubbles}} & \mbox{{\bf relevant product}} & \mbox{{\bf
contributed terms}} \\
\mbox{{\bf in Type II pictures}} & \pi\ast\pi'=(1|2|3)+\cdots &
\mbox{(up to sign)} \\
\hline\hline
a\emp\,\mbox{(three)} & (3|3)\ast(1|0) & 3\{ m(3|3)\}\{\, \{
a\}\{ b,c\} |d,e,f\} \\ 
\emp b,\,\emp c\,\mbox{(two each)} & (3|3)\ast(0|1) & {}\\ \hline
b\emp,\,c\emp\,\mbox{(four each)} & (1|5)\ast(1|0) & 5\{ m(1|5)\}\{ a|
\{ b,c\}\{ d,e,f\}\,\} \\
\emp d,\,\emp e,\,\emp f\,\mbox{(three each)} & (1|5)\ast(0|1) & {} \\ \hline
\end{array}\]
We add up all contributed terms, expanding the ones in the first table
according to the instructions in Section~\ref{ikub}, and set the sum equal
to zero. 
\end{ex}

Examples of lower identities (up to three points in all) can be found in 
\cite{KVZ}.

\begin{rem} The coefficients in Type~II are a major nuisance and are not a 
part of the $\Gi$~identities introduced in~\cite{KVZ}. In order to get rid of
them one may sacrifice consistency in the definition of composition
as a whole and define Type~II as $t-1$ pictures in which an empty 
bubble is skewered to the (say) bottom of any two adjacent lines; then the
two lines protruding above will be treated as if the bubble is full;
i.e. they will be intertwined into one etc. and we will have the 
``non-quadratic'' part without the coefficient. Another way out would be to
define a new unary operation, say~$F$, on the partitioned maps which
describes what happens when adjacent strands are intertwined two by two,
and write the master equation as 
\[ \{\tilde{m}\}\{\tilde{m}\} +F(\tilde{m})=0,\]
where the first summand comes from Type I pictures only,
and $F$ has total degree 1 (mod 2).\label{remar}
\end{rem}

\section{Substructures and examples}

\subsection{Substructures}\label{subs}

Let $(A,m)$ be a $\Gi$ algebra as in Definition~\ref{defiki}. Note that 
if desired, we
may modify this algebraic definition to make any subset of $\{ m(\pi)\}_{
\pi\in\Beb}$, especially of $\{ m(\pi)\}_{\pi\in\Be_0}$, vanish.

\begin{lemma} The infinite sum
\[ m_A=m(0)+m(1)+m(2)+\cdots \]
satisfies the relation
\be \label{AAA} \{\tilde{m}_A\}\{\tilde{m}_A\} =0.\ee
Therefore, $(A,m_A)$ is an $\Ai$ algebra and $(A,l_A)$ is an $\Li$ algebra,
where $l_A$ is the term-by-term graded antisymmetrization of $m_A$.
\end{lemma}

{\it Proof.} The identity (\ref{AAA}) is equivalent to
\[ \{\,\{\tilde{m}\}\{\tilde{m}\}\,\} (i)=0\;\;\;\forall i\geq 0, \]
hence is valid, for the following reasons: the span of $\{ (i)\}_{i\geq 1}$ is
a subalgebra of $(\Pe,\ast)$. Moreover, if
\[ (i)\ra\pi\ast\pi',\]
then we must have
\[ \bar{d}(\pi)+\bar{d}(\pi')=\bar{d}(i)=0,\]
or
\[ \bar{d}(\pi)=\bar{d}(\pi')=0.\Box\] 

\begin{lemma} The infinite sum
\[ m_B=m(1)+m(1|1)+m(1|1|1)+\cdots \]
satisfies
\be \label{BBB} \{\tilde{m}_B\}\{\tilde{m}_B\} =0,\ee
provided that all $\{ m(\pi)\}_{\pi\in\Be_0}$ vanish. Thus $(A,m_B)$ is an
$\Ai$ algebra and $(A,l_B)$ is an $\Li$ algebra, where $l_B$ is the graded
antisymmetrization of $m_B$.
\end{lemma}

{\it Proof.} This time (\ref{BBB}) is equivalent to
\[ \{\,\{\tilde{m}\}\{\tilde{m}\}\,\}\pi_t=0\;\;\;\forall\pi_t=
(1|1|\cdots |1),d(\pi_t)=\bar{d}(\pi_t)=t-1.\]
To see why, we note that the span of $\{\pi_t\}_{t\geq 1}$ is again a
subalgebra of $\Pe$, and whenever
\[ \pi_t\ra\pi\ast\pi',\]
we have
\[ d(\pi)+d(\pi')=\bar{d}(\pi)+\bar{d}(\pi')=t-1\]
and
\[ d(\pi)\geq\bar{d}(\pi),\;\;\; d(\pi')\geq\bar{d}(\pi')\;\;\;
\mbox{(always true for $\pi,\pi'\in\Be$),}\]
implying
\[ d(\pi)=\bar{d}(\pi)\;\;\;{\rm and}\;\;\; d(\pi')=\bar{d}(\pi').\Box\]

In \cite{KVZ}, $(A,m_B)$ in a topological 
vertex operator algebra $A$ (more generally, in a $\Gi$ algebra) is said to
give rise to an $\Li$ algebra in the graded antisymmetrization. In fact,
$m_B$ misses being $\Ai$ itself because of the nonzero maps $m(1|0)$ and
$m(0|1)$: we have
\[ \pi_t\ra\pi\ast\pi'\;\;\;\mbox{with $\pi =(1|1|\cdots |2|\cdots |1|1)$
and $\pi'=(1|0)$ or $(0|1)$,}\]
in addition to the products
\[ \pi_r\ast\pi_s=\pi_{r+s-1}.\]

\subsection{Hochschild complex revisited} \label{hrev}

Recall that the truncated Hochschild complex $\trh$ is a right pre-Lie
algebra under the multiplication
\[ M(1|1)(x,y)\defi\{ x\}\{ y\} ,\]
regardless of any structure on the vector space $A$
(Proposition~\ref{propo2}). Moreover, if $(A,m)$ is an associative algebra,
then the {\bf differential}
\[ M(1)(x)\defi [m,x],\]
{\bf dot product}
\[ M(2)(x,y)\defi\pm\{ m\}\{ x,y\},\]
and partitioned maps
\[ M(1|n)(x,y_1,\dots,y_n)\defi\pm\{ x\}\{ y_1,\dots,y_n\}\;\;\; n\geq 1\]
are known to satisfy the $\Gi$ identities (\cite{GJ2},\cite{KVZ}); this was
first noticed by Gerstenhaber and Voronov and communicated to Getzler. We
provide a brief algebraic proof, in the light of the new terminology.

\begin{prop}[Gerstenhaber-Voronov-Getzler-Jones] The truncated Hochschild 
complex 
\[ (\trh,M(1),M(2),\{ M(1|n)\}_{n\geq 1})\]
is a $\Gi$ algebra.
\end{prop}

{\it Proof.} We have
\bea &&(1)\ast(1)=(1)\nn\\
&&(1)\ast(2)=(2)\ast(1)=(2)\nn\\
&&(2)\ast(2)=(3)\nn\\
&&(1)\ast(1|n)=(1|n)\ast(1)=(1|n)\nn\\
&&(2)\ast(1|n)=(1|n)\ast(2)=(1|n+1)+(2|n)\nn\\
&&(1|n)\ast(1|m)=\mbox{sum of triple-slot terms}\nn\eea
as the only products involving $(1)$, $(2)$, or $(1|n)$ on either side and
leading to a nontrivial identity. But by Remark~\ref{rema} we need only check
\bea &&\mbox{(i)}\;\;\; M(1)^2=0\nn\\
&&\mbox{(ii)}\;\;\; \{ M(1)\}\{ M(2)\}\pm\{ M(2)\}\{ M(1)\} =0\nn\\
&&\mbox{(iii)}\;\;\; (\,\{ M(1)\}\{ M(1|n+1)\}\pm\{ M(1|n+1)\}\{ M(1)\}\nn\\
&&\;\;\;\;\;\;\; \pm\{ M(2)\}\{ M(1|n)\}\pm\{ M(1|n)\}\{ M(2)\}
\, )(1|n+1)=0,\nn\eea
and not
\bea &&\mbox{(iv)}\;\;\;\{ M(2)\}\{ M(2)\} =0,\nn\\
&&\mbox{(v)}\;\;\; (\,\{ M(2)\}\{ M(1|n)\}\pm\{ M(1|n)\}\{ M(2)\}\,
)(2|n)=0,\nn\eea 
nor the triple-slot identities.The first two tell us that the differential
is square-zero and is a derivation of the dot product; these well-known
facts are in \cite{Ger}. The third is Eqn.~(4) of \cite{KVZ} and Eqn.~(8)
of \cite{VG} in disguise. Although not strictly a $\Gi$ identity, Eqn.~(iv)
happens to hold in the complex (the dot product is associative) but
Eqn.~(v) doesn't (it looks like Eqn.~(3) in \cite{KVZ} but contains fewer
terms). Short proofs of (i)-(iv) are also in \cite{A}. $\Box$

\begin{rem} Once again, the result holds for the full complex with all
$\{ a\}\{ b_1,\dots,b_n\} =0$ for $a$, $b_i\in A$.
It would also be interesting to study the cases where $A$ is an
$\Ai$ algebra (\cite{Get},\cite{A}) and/or $\trh$ is replaced by $\CAA$
with the additional composition maps $M(\la)$ \cite{ABV}.
\end{rem}

\subsection{Phi-operators revisited} \label{phiope}

The multilinear operators $\Phi_{\De}^r$ defined at the end of
Section~\ref{ikub} may be a convenient tool to construct brand new $\Gi$
algebras starting with a small data set. Let us consider an easy truncated
version: define
\bea &&\Phi(1)\defi\Phi_{\De}^1=\De\nn\\
&&\Phi(2)=\Phi(1|1)\defi\Phi_{\De}^2\nn\\
&&\Phi(1|2)\defi\Phi_{\De}^3,\nn\eea
where $\De$ is any odd square-zero operator on an algebra $(A,m)$ ($m$ not
necessarily commutative or associative). This truncation especially makes
sense when $\De$ is a differential operator of order three. The relevant
products in $\Pe$ are
\bea &&(1)\ast(1)=(1)\nn\\
&&(1)\ast(2)=(2)\ast(1)=(2)\nn\\
&&(1)\ast(1|1)=(1|1)\ast(1)=(1|1)\nn\\
&&(1)\ast(1|2)=(1|2)\ast(1)=(1|2)\nn\\
&&(2)\ast(1|1)=(1|1)\ast(2)=(1|2)+(2|1),\nn\eea 
and we want
\bea &&\mbox{(i)}\;\;\;\De^2=0\nn\\
&&\mbox{(ii)}\;\;\;\mbox{$\De$ is a derivation of $\Phi_{\De}^2$}\nn\\
&&\mbox{(iii)}\;\;\;\mbox{ditto}\nn\\
&&\mbox{(iv)}\;\;\; (\,\{\Phi(1)\}\{\Phi(1|2)\}\pm\{\Phi(1|2)\}\{\Phi(1)\}
\nn\\ &&\;\;\;\;\;\;\;\;\;\pm
\{\Phi(2)\}\{\Phi(1|1)\}\pm\{\Phi(1|1)\}\{\Phi(2)\}\, )(1|2)=0.\nn\eea
Again (i)-(iii) are known (see \cite{ABV} or \cite{A}), and Eqn.~(iv) can
be obtained from Lemma~6 of \cite{A}:

\begin{lemma} For odd linear operators $T$ and $U$ on $A$, the Gerstenhaber
bracket ${[T,U]}=TU+UT$ is related to the Gerstenhaber bracket of the
$\Phi$ operators as follows:
\bea \Phi_{[T,U]}^1(a)&=&[\Phi_T^1,\Phi_U^1]\{ a\}\nn\\
\Phi_{[T,U]}^2(a,b)&=&[\Phi_T^1,\Phi_U^2]\{ a,b\} +[\Phi_U^1,\Phi_T^2]\{
a,b\} \nn\\
\Phi_{[T,U]}^3(a,b,c)&=&[\Phi_T^1,\Phi_U^3]\{ a,b,c\} +[\Phi_U^1,\Phi_T^3]
\{ a,b,c\}\nn\\
&&+[\Phi_T^2,{\rm ad}(\Phi_U^2)\{ a\} ]\{ b,c\} +[\Phi_U^2,{\rm ad}(\Phi_T^2
)\{ a\} ]\{ b,c\} .\nn\eea
Here the adjoint operator is defined as
\[ {\rm ad}(\Phi_{\De}^r)\{ a_1\}\,\{ a_2,\dots,a_r\}\defi\Phi_{\De}^r(a_1,
a_2,\dots,a_r).\]
\end{lemma}

By substituting $T=U=\De$ we find that the last equation in the above Lemma
is exactly Eqn.~(iv), by virtue of ${[\De,\De]}=2\De^2=0$. Therefore, we
have

\begin{prop} If $\De$ is an odd square-zero operator on an arbitrary
algebra $(A,m)$, then
\[ (A,\Phi(1),\Phi(2),\Phi(1|1),\Phi(1|2))\]
is a $\Gi$ algebra.
\end{prop}

\subsection{Topological vertex operator algebras}

Topological vertex operator algebras (TVOA) have always fueled the subject of
$\Gi$ algebras; see \cite{LZ} and \cite{KVZ}. The article \cite{KVZ}
indicates the existence of a $\Gi$ structure on a TVOA but does not produce
the actual products. 
We will assume some acquaintance with vertex operator super algebras (VOSA)
and TVOA's, which are VOSA's with
extra structure (\cite{FLM}, \cite{FHL}, \cite{DL}, 
\cite{Geb}, \cite{LZ2}, \cite{LZ3}, \cite{LZ4}). A {\bf 
VOSA} is a {\bf Z}-bigraded vector space~$V$ (one $L_0$ and one super, or
fermionic, grading) in which we associate to every element $u$ a unique
{\bf vertex operator} $u(z)=\sum u_nz^{-n-1}$, with $u_n\in End(V)$. There
is an action of the Virasoro algebra by $L(z)=\sum L_nz^{-n-2}$ where the
eigenvalues of $L_0$ are bounded from below and $L_{-1}$ is formal
differentiation. The vacuum element {\bf 1}, represented by $1\cdot z^0$,
is the identity element with respect to the distinguished {\bf Wick
product} given by $u_{-1}v$. Although the overall identities satisfied by
all the bilinear products $u_nv$ are summed up by the formal relation ([DL])
\[ {[u(z_1),v(z_2)]}(z_1-z_2)^t=0\;\;\;\mbox{for sufficiently large 
$t=t(u,v)$,}\]
the specialized identities
\[ (u_mv)_n=\sum_{i\geq 0}(-1)^i\mi(u_{m-i}v_{n+i}-(-1)^{m+|u||v|}
v_{m+n-i}u_i)\]
and
\[ {[u_m,v_n]}=\sum_{i\geq 0}\mi(u_iv)_{m+n-i}\]
(for $m$, $n\in${\bf Z}) are very useful. For the record, we list a number
of products already shown to satisfy the $\Gi$ identities:
\begin{itemize}
\item The differential (``BRST operator''), usually denoted by $Q$,
can be taken to be the odd linear operator $m(1)$.
\item The Wick product (``normal ordered product'') plays the role of the
even bilinear operator $m(2)$ in any VOSA. 
\item The odd trilinear operation $n(u,v,t)$ of \cite{LZ} (Eqn.~(2.16))
looks like $m(3)$.
\item The odd bilinear product $m(u,v)$ of \cite{LZ} (Eqn.~(2.14)) is like 
$m(1|1)$.
\end{itemize}
There are additional products in \cite{LZ} and \cite{ABV} on a TVOA which
may eventually be linked to the $\Gi$ structure (most of these are related
to homotopy Batalin-Vilkovisky algebras).
\vsp

In addition, we establish a surprising property of the Wick product on a
VOSA:
 
\begin{prop} Any VOSA $V$ with the Wick product is a left pre-Lie
algebra. Therefore, the {\bf Wick bracket} defined by
\[ {[u,v]}_W\defi u_{-1}v-(-1)^{|u||v|}v_{-1}u\]
is a graded Lie bracket on $V$.
\end{prop}

{\it Proof.} Dropping the super degrees, we have
\[ (u_{-1}v)_{-1}=\sum_{i\geq 0}(u_{-1-i}v_{-1+i}+v_{-2-i}u_i)=u_{-1}
v_{-1}+\sum_{i\geq 0}(u_{-2-i}v_i+v_{-2-i}u_i),\]
so that for the left (Wick) multiplication operators $L_u$ we have the 
identity
\bea &&L_uL_v-L_vL_u-L_{[u,v]_W}\nn\\
&=&u_{-1}v_{-1}-v_{-1}u_{-1}-(u_{-1}v)_{-1}+(v_{-1}u)_{-1}\nn\\
&=&u_{-1}v_{-1}-v_{-1}u_{-1}-u_{-1}v_{-1}-\sum_{i\geq 0}(u_{-2-i}v_i+
v_{-2-i}u_i)+v_{-1}u_{-1}+\sum_{i\geq 0}(v_{-2-i}u_i+u_{-2-i}v_i)\nn\\
&=&0.\;\;\;\Box\nn\eea
\vsp

It was pointed out to the author by Haisheng Li that the Jacobi identity
for ${[u,v]}_W$ was independently observed by Dong, Li, and Mason in
\cite{DLM}. The author has also found some evidence (in low $L_0$ degrees) 
that the {\bf
Zhu product}~$\ast_0$ \cite{Z} on a VOSA will turn out to be a left pre-Lie
product; a direct proof seems formidable (meanwhile, higher products
$\ast_n$ for $n\geq 1$ as defined by Dong, Li, and Mason in \cite{DLM2} do
not seem to be pre-Lie). The Zhu product makes a quotient of any VOSA into
an associative algebra, and it is interesting that the {\bf Zhu bracket}
\[ {[u,v]}_Z\defi u\ast_0 v-(-1)^{|u||v|}v\ast_0 u\]
may already be Lie on the original space. The pre-Lie identity is also part
of the definition of a {\bf Novikov algebra} (see e.g. Osborn's \cite{Os})
which has been studied in some detail in terms of the representation theory
and also in the context of vertex operator algebras, variational calculus,
etc. 
\vsp

Finally, since the modes $u_n$ of vertex operators for $n\geq 0$ have been
identified as differential operators of order $n+1$ with respect to the
Wick product, and the $\Phi$ operators have been shown to be an invaluable
tool in proving identities on VOSA's (such as those related to the
generalized Batalin-Vilkovisky structure) (\cite{ABV},\cite{A}), we expect
the $\Phi$ operators to play some role in making explicit the $\Gi$
structure on a TVOA.


\end{document}